\documentclass[12pt,a4paper, twoside]{article}

\usepackage{amsmath,amssymb,amsthm,amsfonts,mathrsfs,amscd,environ}
\usepackage{dsfont}
\usepackage{latexsym,enumerate,color,geometry,extarrows}
\usepackage{xcolor}
\usepackage{verbatim,fancyhdr,enumitem}

 \NewEnviron{ews}{%
\begin{equation}\begin{split}
  \BODY
\end{split}\end{equation}
}

\NewEnviron{ews*}{%
\begin{equation*}\begin{split}
  \BODY
\end{split}\end{equation*}
}

\def\beg{\begin}
\def\bequ{\begin{equation}}
\def\enqu{\end{equation}}
\def\bes{\begin{split}}
\def\ens{\end{split}}
\def\bews{\begin{ews}}
\def\beqn{\begin{eqnarray}}
\def\enqn{\end{eqnarray}}
\def\beq*{\begin{equation*}}
\def\enq*{\end{equation*}}
\def\bqn*{\begin{eqnarray*}}
\def\eqn*{\end{eqnarray*}}
\def\bary{\begin{array}}
\def\eary{\end{array}}
\def\bpma{\begin{pmatrix}}
\def\epma{\end{pmatrix}}
\def\bvma{\begin{Vmatrix}}
\def\evma{\end{Vmatrix}}

 \numberwithin{equation}{section}

\def\al{\alpha}
\def\be{\beta}
\def\ga{\gamma}
\def\de{\delta}
\def\ep{\epsilon}
\def\vep{\varepsilon}

\def\la{\lambda}

\def\si{\sigma}
\def\ta{\tau}
\def\ph{\phi}

\def\Om{\Omega}

\def\Q{\mathbb Q}
\def\R{\mathbb R}
\def\P{\mathbb P}
\def\E{\mathbb E}

\def\sF{\mathscr F}

\def\e{\operatorname{e}}

\def\d{\mathrm{d}}

\def\ff{\frac}
\def\ra{\rightarrow}

\def\<{\langle}
\def\>{\rangle}
\def\sq{\sqrt}

\def\we{\wedge}
\def\1{\mathds{1}}

\def\8{\infty}

\allowdisplaybreaks

\title{{\bf Weak convergence of Euler scheme for SDEs with singular drift}
}

\author{
{\bf Yongqiang Suo$^{a)}$, Chenggui Yuan$^{a)}$, Shao-Qin Zhang$^{b)}$}\\
\footnotesize{$^{a)}$Mathematic department, Swansea University, Bay campus, SA1 8EN, UK}\\
\footnotesize{Email: 971001@swansea.ac.uk\qquad C.Yuan@swansea.ac.uk}\\
\footnotesize{$^{b)}$School of Statistics and Mathematics}\\
 \footnotesize{Central University of Finance and Economics, Beijing 100081, China}\\
\footnotesize{Email: zhangsq@cufe.edu.cn}\\
}

\begin{document}

\maketitle

\begin{abstract}
In this paper, we investigate the weak convergence rate of Euler-Maruyama's  approximation for stochastic differential equations with irregular drifts. Explicit weak convergence rates are presented if  drifts satisfy an integrability condition including discontinuous functions which can be non-piecewise continuous or in fractional Sobolev space.
\end{abstract}\noindent

\noindent


\section{Introduction}
Stochastic differential equations (SDEs for short) with singular coefficients have been extensively studied recently, see \cite{NM,W18,XZ16,XZ,X1} and references therein. Meanwhile, in order for one to understand the numerical approximation of SDEs with irregular coefficients, numerical schemes are established.   The strong and weak convergence rate for singular SDEs are obtained,  see \cite{BHY,BHZ,KLY,SV,GM1,GM2,LM,AML,NT1,NT2,S} for instance. Particularly, the strong convergence of  Euler-Maruyama's (abbreviated as EM's) scheme for discontinuous monotone drifts was investigated in \cite{NP}. \cite{NT1} obtained the strong convergence rates of EM's scheme for SDEs with drifts satisfying the one-side Lipschitz condition,  and \cite{NT2} investigated the one-dimensional setup without the one-side Lipschitz assumption. Recently, \cite{BHZ} obtains strong convergence rates for multidimensional SDEs under an integrability condition by using the Krylov's estimate and the result of Gaussian type heat kernel estimates established by the parametrix method in \cite{LM}.  The weak convergence is concerned with the convergence of the distribution of the solutions of SDEs.  In this paper, we shall investigate the weak error of  EM's scheme for the following SDE on $\R^d$ 
\beg{align}\label{eq2.1}
\d X_t=  b(X_t)\d t+\si \d W_t,~X_0=x\in\R^d,
\end{align}
where $(W_t)_{t\ge0}$ is a Brownian motion on $\R^d$ with respect to a complete filtration probability space $(\Om, (\sF_t)_{t\ge0}, \sF, \P)$. The associated EM's scheme reads as follows: for any $\de \in (0, 1),$
\beg{align}\label{eq2.2}
\d X_t^{(\de)}=b(X_{t_\de}^{(\de)})\d t+\si\d W_t,~X_0^{(\de)}=0,
\end{align}
where $t_\de=[t/\de]\de$ and $[t/\de]$ denotes the integral part of $t/\de$.  The weak convergence rate is concerned with the approximation of $\E f(X_t)$ by $\E f(X_t^{(\de)})$ for a given function $f$. The weak error has been obtain for some SDEs with discontinuous drifts in \cite{KLY,SV}. However, they required that the given function $f$ is H\"older continuous. The weak convergence rate with a measurable and bounded function $f$  can be dated back to \cite{BT}, where the coefficients of SDEs need to be smooth. Recently, \cite{B,S} established the weak convergence rate of EM's scheme for SDEs with irregular coefficients by using Girsanov's transformation. Inspired by \cite{BHZ} and \cite{B,S}, we shall give a  note on the weak error for \eqref{eq2.3} with a possibly discontinuous drift $b$. Moreover, the given function $f$ is only assumed to be bounded and measurable on $\R^d$.

The remainder of this paper is organized as follows: The main result is presented in Section 2. All the proofs are given in Section 3.

\section{Main Result and Examples}

Let $|\cdot|$ be the Euclidean norm, $\langle\cdot, \cdot \rangle$ be the Euclidean product.  $\|\cdot\|$ denotes the operator norm.  Throughout this paper, we assume the coefficients of \eqref{eq2.3} satisfy the following assumptions:
\beg{enumerate}
\item[(H1)] $b:\R^d\rightarrow\R^d$ is measurable and $\si$ is an invertible $d\times d$-matrix. There exist nonnegative constants $L_1,L_2 $ such that 
$$|b(x)|\leq L_1+L_2|x|.$$

\item[(H2)] There exist $p_0\geq 2$, $\al>0$  and $\ph\in C([0,T];(0,+\infty))$ with $\int_{0^+}\ph^2(s)\d s<\infty$ such that
\beg{align*}
\sup_{z\in\R^d}\int_{\R^d\times\R^d}|b(y)-b(x)|^{p_0}\ff {\e^{-\ff {|x-z|^2} s-\ff {|y-x|^2} r}} {s^{\ff d 2}r^{\ff d 2}}\d x\d y\leq (\ph(s)r^{\al})^{p_0} ,~s>0,r\in [0,1].
\end{align*}
\end{enumerate}

By  \cite[ Theorem 1.1]{X1}, \eqref {eq2.1} has a unique strong solution under (H1).  It is clear that \eqref{eq2.2} also has a unique strong solution. We denote  $\|f\|_{\infty}=\sup_{x\in\R^d}|f(x)|$. We now formulate the main result.
\beg{thm}\label{2}
Assume (H1)-(H2). If 
\beg{align}\label{ine-Tp0}
 TL_2 \|\si^{-1}\| \|\si\| \ff {\sq{2(p_0+1)(p_0+3)}} {p_0-1}<1,
\end{align}
then for any bounded measurable function $f$ on $\R^d$, there exists a constant $C_{T,p_0,\si,x}>0$ such that  
\beg{align}\label{eq2.3}
|\E f(X_t)-\E f(X_t^{(\de)})|\le C_{T,p_0,\si,x}\|f\|_\infty\de^{\al},  t\in[0,T].
\end{align}
Additionally,  if $b$ has sublinear growth, i.e.  for any $\ep>0$, there exists $L(\ep)>0$ such that
$|b(x)|\leq L(\ep)+\ep|x|$, the convergence holds for any $T>0$.
\end{thm}

\beg{rem} {\rm When the drift $b$ is non-regular, the boundedness on $b$ is needed see e.g. \cite{BHZ,KLY,SV}. 
 Here, we allow that $b$ has linear growth by (H1).  If $b$ is bounded or $\be$-H\"older continuous with $\be<1$, then $b$ has sublinear growth. 

In the condition (H2), if $\al$ is a decreasing function of  $p_0$,  then we can choose  $p_0=2$ without considering that $T$ depends on $p_0$ increasingly, see  Example \ref{e2}. 
}
\end{rem}
We give several examples to illustrate  the condition (H2) and the convergence rate $\al$.

\beg{exa}\label{e1}
If $b$ is the H\"older continuous with exponent $\be$, i.e. 
$$|b(y)-b(x)|\le L|x-y|^\be,$$ 
then (H2) holds with $\al=\ff{\be}{2}$ and  a constant function $\phi(s).$  It is clear that $b$ has sublinear growth if $\be<1$. Then for any $T>0$, \eqref{eq2.3} holds with $\al=\ff {\be } 2$. 
\end{exa}
\beg{proof} By the H\"older continuity and the fact 
\beg{align}\label{ineq}
\sup_{x\ge0}(x^{\ga'}\e^{-\ga x^2})=\Big(\ff{{\ga'}}{2\e\ga}\Big)^{{\ga'}/2},~~~~~~{\ga'},\ga>0,
\end{align}
the assertion follows from the following inequality
\beg{align*}
&\sup_{z\in\R^d}\int_{\R^d\times\R^d}|b(y)-b(x)|^p\ff {\e^{-\ff {|x-z|^2} s-\ff {|y-x|^2} r}} {s^{\ff d 2}r^{\ff d 2}}\d x\d y\\
&\le L^p\sup_{z\in\R^d}\int_{\R^d\times\R^d}|y-x|^{\be p}\ff {\e^{-\ff {|x-z|^2} s-\ff {|y-x|^2} r}} {s^{\ff d 2}r^{\ff d 2}}\d x\d y\\
&\le L^p\ff{1}{{s^{\ff d 2}r^{\ff d 2}}}\bigg(\ff{\be pr}{\e}\bigg)^{\ff{\be p}{2}}\sup_{z\in\R^d}\int_{\R^d\times\R^d}\e^{-\ff {|x-z|^2}{s}}\e^{-\ff {|y-x|^2}{2r}}\d x\d y\\
&\le CL^p \left(\ff{\be pr}{\e}\right)^{\ff{\be p}{2}}.
\end{align*}

\end{proof}

The following example shows that (H2) can hold even if $b(x)$ is not piecewise continuous.
\beg{exa}\label{e2}
 Let $A$ be the Smith-Volterra-Cantor set on $[0,1]$, which is constructed in the following way. The first step, we let $I_{1,1}=\left(\ff 3 8,\ff 5 8\right)$, $J_{1,1}=\left[0,\ff 3 8\right],$ $J_{1,2}=[\ff 5 8,1]$ and  remove the open interval $I_{1,1}$. The second step, we  remove the middle $\ff 1 {4^2}$ open intervals, denoting by $I_{2,1}$ and $I_{2,2}$, from  $J_{1,1}$ and $J_{1,2}$ respectively, i.e.  $I_{2,1}=\left(\ff 5 {32}, \ff {7} {32}\right)$, $I_{2,2}=\left(\ff {25} {32},\ff {27} {32}\right)$. The intervals left are denoted by $J_{2,1},J_{2,2},J_{2,3},J_{2,4}$, i.e.
$$J_{2,1}=\left[0,\ff 5 {32}\right],J_{2,2}=\left[\ff 7 {32},\ff 3 8\right],J_{2,3}=\left[\ff 5 {8}, \ff {25} {27}\right],J_{2,4}=\left[\ff {27} {32},1\right].$$ 
For the $n$-th step, we  remove the middle $\ff 1 {4^{n}}$ open intervals $I_{n,1},\cdots, I_{n,2^{n-1}}$ from  $J_{n-1,1},\cdots, J_{n-1,2^{n-1}}$ respectively, and the intervals left are denoted by $J_{n,1},\cdots, J_{n, 2^n}$. Let
$$A=\bigcap_{n=1}^\infty\left( \bigcup_{k=1}^{2^n} J_{n,k}\right).$$
Then $A$ is a nowhere dense set and the Lebesgue measure of $A$ is $1/2$. Define
\beg{align*}
b(x)&=\1_{[0,1]}(x)-\sum_{n=1}^\infty \sum_{j=1}^{2^{n-1}}2^{-(n+j)}\1_{I_{n,j}}(x)\\
&=\1_A(x)+\sum_{n=1}^\infty \sum_{j=1}^{2^{n-1}}\left(1-2^{-(n+j)}\right)\1_{I_{n,j}}(x).
\end{align*}
All of the endpoints of the intervals $\bar I_{n,j}$ are the discontinuous points of $b$, which is dense in $A$.  For any  interval $I\subset [0,1]$ such that $I\cap A\neq \emptyset$, it always contains the discontinuous points of $b$. However, any interval $I\subset [0,1]$ such that $I\cap A=\emptyset$,  it is a subset of some $I_{n,j}$. Hence, $b$ is not a piecewise continuous function. In the following, we shall show that $b$ satisfies condition (H2)  with $p_0=2$ and $\al=\ff 1 4$ and $\ph(s)=Cs^{\ff 1 2}$.
\end{exa}
\beg{proof}
For $z>0$ and any interval $(a_1,a_2)$ (it is similar for $[a_1,a_2]$)
\beg{align*}
&\int_{-\infty}^{+\infty}  \left|\1_{(a_1,a_2)}(x+u)-\1_{(a_1,a_2)}(x)\right|^2 \d x\\
&\qquad =\int_{a_1-u}^{a_2-u}\1_{(a_1,a_2)^c}(x)\d x+\int_{a_1}^{a_2}\1_{(a_1-u,a_2-u)^c}(x)\d x\\
&\qquad =\int_{a_1-u}^{(a_2-u)\wedge a_1}\d x+\int_{(a_2-u)\vee a_1}^{a_2}\d x\\
&\qquad  \le 2\left( |u|\we (a_2-a_1) \right).
\end{align*}
For $z<0$, 
\beg{align*}
&\int_{-\infty}^{+\infty} \left|\1_{(a_1,a_2)}(x+u)-\1_{(a_1,a_2)}(x)\right|^2 \d x\\
&\qquad = \int_{-\infty}^{+\infty} \left|\1_{(a_1,a_2)}(v)-\1_{(a_1,a_2)}(v-u)\right|^2 \d v\le 2\left( |u|\we (a_2-a_1) \right).
\end{align*}
Hence, by Jessen's inequality
\beg{align*}
& \int_{-\infty}^{+\infty}|b(x+u)-b(x)|^2\d x\\
&\qquad \leq  \int_{-\infty}^{+\infty}\left(\left|\1_{[0,1]}(x+u)-\1_{[0,1]}(x)\right|\right.\\
&\qquad\qquad\left. +\sum_{n=1}^{+\infty}\sum_{j=1}^{2^{n-1}}2^{-(n+j)}\left|\1_{I_{n,j}}(x+z)-\1_{I_{n,j}}(x)\right|\right)^2\d x\\
&\qquad\leq \left(1+\sum_{n=1}^{+\infty}\sum_{j=1}^{2^{n-1}}2^{-(n+j)}\right)\left\{ \int_{-\infty}^{+\infty} \left|\1_{[0,1]}(x+u)-\1_{[0,1]}(x)\right|^2\d x\right.\\
&\qquad\qquad \left. +\sum_{n=1}^{+\infty}\sum_{j=1}^{2^{n-1}}2^{-(n+j)}\int_{-\infty}^{+\infty}  \left|\1_{I_{n,j}}(x+u)-\1_{I_{n,j}}(x)\right|^2 \d x\right\}\\
& \qquad \leq 2\left(1+\sum_{n=1}^\infty\sum_{j=1}^{2^{n-1}}2^{-(n+j)}\right)^2 |u|
 =4|u|
\end{align*}
Therefore
\beg{align*}
&\sup_{z\in\R}\int_{\R\times\R}|b(y)-b(x)|^2\frac{{\e^{-\ff {|x-z|^2} s}\e^{-\ff {|y-x|^2} r}}} {s^{\ff 1 2}r^{\ff 1 2}}\d x\d y\\
&\qquad\le\ff{1}{{s^{\ff 1 2}r^{\ff 1 2}}}\int_{\R}\e^{-\ff {|u|^2} r}\int_{\R}|b(x+u)-b(x)|^2\d x\d u\\
&\qquad\le\ff{4}{{s^{\ff 1 2}r^{\ff 1 2}}}\int_{\R}\e^{-\ff {|u|^2} r}|u|\d u
=\left(Cs^{-\ff 1 4}r^{\ff 1 4}\right)^2.
\end{align*}
\end{proof}

A general class of functions that satisfies (H2) is the (fractional) Sobolev space $W^{\be,p}(\R^d)$:

\beg{exa}\label{e3}
If there exist $\be >0$ and $p\in (2,\infty)\cap (\ff d 2,+\infty)$ such that the Gagliardo seminorm of $b$ is finite, i.e.
\beg{align*}
\left[b\right]_{W^{\be,p}}:=\left(\int_{\R^d\times\R^d} \ff {|b(x)-b(y)|^p} {|x-y|^{d+\be p}}\d x\d y\right)^{\ff 1 p}<\infty,
\end{align*}
then (H2) holds for any $2\leq p'<p$ with $\al=\ff {\be} 2$ and $\ph(s)= Cs^{\ff {d} {2p}}$.  Hence, if $b$ satisfies (H1) and $\left[f\right]_{W^{\be,p}}<\infty$ with $p\in (2,\infty)\cap (\ff d 2,+\infty)$, then \eqref{eq2.3} holds.
\end{exa}
The proof of this example is similar to that of example in \cite{BHZ}, we omit it here.

\beg{rem} {\rm In \cite{BHZ}, the strong convergence and the convergence rate are investigated under the drift satisfying an integrability condition and boundeness. Here we obtain  the weak convergence rate, which the drift needs not be bounded and the function $f$ in equation \eqref{eq2.3} is only bounded and measurable. From the examples above, one can see that the drift could be very irregular, this means that we have extended the results   in \cite{BT} where the coefficients must be smooth.
 However, the  method developed here seems  difficult to  deal with SDEs driven by multiplicative noise. Moreover, our method is not optimal in Lipschitz case since 
  the classical weak rate is $\al=1$ for SDEs with smooth coefficients in \cite{BT}. }
 \end{rem}

\section{Proof of Theorem \ref{2}}
The key point to prove the main result is to construct a reference SDE, which provides a new representation of \eqref{eq2.1} and  its EM's approximation SDE \eqref{eq2.2} under another probability measures which will be defined by the Girsanov theorem.  

Let $Y_t=\si W_t$, which is the reference SDE of \eqref{eq2.1}.  The first lemma is on the exponential estimate of $|b(Y_t)|$. We use a weaker condition than that of (H1) on $b$:  there exist  nonnegative constants $L_1,L_2 $ and  $F\geq 0$ satisfying $\|F\|_{L^{p_1}(\R^d)}$ for some $p_1> d$ such that 
\beg{align}\label{ine-bn}
|b(x)|\leq L_1+L_2|x|+F(x).
\end{align}
By Krylov's estimate, for any $q$ such that $\ff d {p_1} +\ff 1 q<1$, there exists (see e.g. \cite{NM})
$$\E \left[\int_S^T F^2(Y_s)\d s\Big| \sF_S\right]\leq (T-S)^{\ff 1 q}\|F\|_{L^{p_1}},$$
which yields the following Khasminskii's estimate (see e.g. \cite[Lemma 3.5]{XZ}): for any $C>0$ 
\beg{align}\label{Khas}
\E\exp\left\{C\int_0^T F^2(Y_s)\d s\right\}<\infty.
\end{align}

\beg{lem}\label{lemn}
Assume \eqref{ine-bn} holds. If $\la$ and $T$ satisfy 
\begin{align}\label{ine-Tla}
2T^2\la L_2^2\|\si^{-1}\|^2\|\si\|^2<1,
\end{align}
then
\beg{align}\label{eqnm}
\E\exp\Big\{\la\int_0^{T_0}|\si^{-1}b(Y_s)|^2\d s\Big\}<\8.
\end{align}
\end{lem}
\beg{proof}
By the H\"older inequality, we derive that
\beg{align}\label{ine-exp-Y}
&\E\exp\Big\{\la\int_0^{T}|\si^{-1}b(Y_s)|^2\d s\Big\}\nonumber\\
&\le \E\exp\Big\{\la \int_0^{T} \|\si^{-1}\|^2 \left((L_1+L_2|x|)+L_2|Y_s-x|+F(Y_s)\right)^2\d s\Big\}\nonumber\\
&\le \exp\{\la T \|\si^{-1}\|^2(L_1+L_2|x|)^2\left(2+\vep_1^{-1}\right)\}\nonumber\\
&\qquad \times\left(\E\exp\left\{\la (1+\vep_1+\vep_2)^2 L_2^2\|\si^{-1}\|^2\int_0^T|Y_s-x|^2\d s\right\}\right)^{  \ff{1}{1+\vep_1+\vep_2}}\nonumber\\
&\qquad \times \left(\E\exp\left\{ \ff{\la(2+\vep_2^{-1})(1+\vep_1+\vep_2)}{\vep_1+\vep_2}\|\si^{-1}\|^2\int_0^T F^2(Y_s)\d s\right\}\right)^{\ff {\vep_1+\vep_2}{1+\vep_1+\vep_2}}.
\end{align}

It follows from \eqref{Khas} that for any $\vep_2>0$
$$\E\exp\left\{\ff{2\la(1+\vep_2^{-1})(1+\vep_1+\vep_2)}{\vep_1+\vep_2}\|\si^{-1}\|^2\int_0^T F^2(Y_s)\d s\right\}<\infty.$$
Since \eqref{ine-Tla}, we can choose $\vep_1$ and $\vep_2$ such that 
$$2T^2(1+\vep_1+\vep_2)^2\la L_2^2\|\si^{-1}\|^2\|\si\|^2<1.$$
Then, by the Jenssen inequality,
\beg{align}\label{ine-Y-x}
&\E\exp\left\{\la (1+\vep_1+\vep_2)^2 L_2^2\|\si^{-1}\|^2\int_0^T|Y_s-x|^2\d s\right\}\nonumber\\
&\qquad \leq \ff 1 T\int_0^T\E\exp\left\{T\la (1+\vep_1+\vep_2)^2 L_2^2\|\si^{-1}\|^2 |Y_s-x|^2 \right\}\d s\nonumber\\
&\qquad = \int_0^T\int_{\R^d}\ff {\exp\left\{T\la (1+\vep_1+\vep_2)^2L_2^2\|\si^{-1}\|^2 |y-x|^2 -\ff {|\si^{-1}(y-x)|^2} {2s}\right\}} {T\sq{(2s\pi)^{d}\det(\si\si^*)}}\d y \d s\nonumber\\
&\qquad <\infty.
\end{align}
The proof is therefore complete.
\end{proof}

For the process $\{Y_{t_\de}\}_{t\in [0,T]}$, we have the following estimate.

\beg{lem}\label{lemm}
Assume (H1) holds.  If $\la$ and $T$ satisfy \eqref{ine-Tla}, then 
\beg{align}\label{eqnn}
\sup_{\de >0}\E\exp\Big\{\la\int_0^{T}|\si^{-1}b(Y_{s_\de})|^2\d s\Big\}<\8.
\end{align}
\end{lem}
\beg{proof}
In the same way as in \eqref{ine-exp-Y}, we have for any $\vep_1>0$,
\beg{align*}
\E\exp\Big\{\la\int_0^{T}|\si^{-1}b(Y_{s_\de})|^2\d s\Big\} &\le \exp\{\la T \|\si^{-1}\|^2(L_1+L_2|x|)^2\left(1+\vep_1^{-1}\right)\}\nonumber\\
&\quad \times\E\exp\left\{\la (1+\vep_1) L_2^2\|\si^{-1}\|^2\int_0^T|Y_{s_\de}-x|^2\d s\right\}.
\end{align*}
By the Jenssen inequality, as \eqref{ine-Y-x}, we can  choose $\vep_1$ such that 
$$\sup_{\de>0}\E\exp\left\{\la (1+\vep_1) L_2^2\|\si^{-1}\|^2\int_0^T|Y_{s_\de}-x|^2\d s\right\}<\infty.$$
\end{proof}
 
We can now  give the Proof of Theorem \ref{2}. Let 
\beg{align*}
\hat W_t &= W_t-\int_0^t\si^{-1}b(Y_s)\d s,\quad \tilde{W}_t  =W_t-\int_0^t\si^{-1} b(Y_{s_\de}) \d s,\\
R_{1,T}&=\exp\Big\{\int_0^T\<\si^{-1}b(Y_s),\d W_s\>-\ff{1}{2}\int_0^T|\si^{-1}b(Y_s)|^2\d s\Big\},\\
R_{2,T}&=\exp\Big\{\int_0^T\<\si^{-1} b(Y_{s_\de}),\d W_s\>-\ff{1}{2}\int_0^T\big|\si^{-1} b(Y_{s_\de})\big|^2\d s\Big\}.
\end{align*}

\beg{proof}[Proof of Theorem \ref{2}] The proof is divided into two steps:

Step (i), we shall prove that the assertion holds under (H1) and (H2).

We first show that $\{\hat W_t\}_{t\in [0,T]}$ is a Brownian motion under $\Q_1:=R_{1,T}\P$, and $\{\tilde W_t\}_{t\in [0,T]}$ is a Brownian motion under $R_{2,T}\P$. Since \eqref{ine-Tp0}, we have by Lemma \ref{lemn} that 
\begin{align}\label{ine-exp}
\E\exp\left\{\ff {(p_0+3)(p_0+1)} {(p_0-1)^2}\int_0^{T}|\si^{-1}b(Y_s)|^2\d s\right\}<\infty.
\end{align}
It is clear that $\ff {(p_0+3)(p_0+1)} {(p_0-1)^2}>\ff 1 2$, so by Novikov's condition $\{R_{1,t}\}_{t\in [0,T]}$ is a martingale and $\{\hat W_t\}_{t\in [0,T]}$ is a Brownian motion under $\Q_1$.  Similarly, it follows from \eqref{ine-Tp0}, Lemma \ref{lemm} and Novikov's condition that $\{\tilde W_t\}_{t\in [0,T]}$ is a Brownian motion under $\Q_2:=R_{2,T}\P$.

We can reformulate $Y_t$ as follows:
$$Y_t=x+\int_0^t b(Y_s)\d s+\si \hat W_t,$$
which means that $Y_t$ is a weak solution of \eqref{eq2.1}. Hence,  $Y_t$ under $\Q_1$ has the same law of $X_t$ under $\P$, since the pathwise uniqueness of the solutions to \eqref{eq2.1}.  Similarly,  reformulating $Y_t$  as follows:
 \beg{align}\label{eq2.5}
Y_t= x+\int_0^t b(Y_{s_\de}) \d s+\si \tilde{W}_t,
 \end{align}
$(Y_t,\tilde{W}_t)$ under $\Q_2$ is also a solution of \eqref{eq2.2}, which has a pathwise unique solution. Hence  $Y_t$ under $\Q_2$ has the same law of $X_t^{(\de)}$ under $\P$.

For every bounded measurable function $f$ on $\R^d$
\beg{align}\label{fr}
&|\E f(X_t)-\E f(X_t^{(\de)})|=|\E_{\Q_1}f(Y_t)-\E_{\Q_2}f(Y_t)|\nonumber\\
&=\E|(R_{1,T}-R_{2,T})f(Y_t)| \le\|f\|_\8\E|R_{1,T}-R_{2,T}|\nonumber\\
&\le\|f\|_\8\E\Big\{(R_{1,T}\vee R_{2,T})\Big|\int_0^T\<\si^{-1}(b(Y_s)-b(Y_{s_\de})),\d W_s\>\nonumber\\
&\qquad\qquad +\ff{1}{2}\int_0^T\Big(|\si^{-1}b(Y_{s_\de})|^2-|\si^{-1}b(Y_{s})|^2\Big)\d s\Big|\Big\}\nonumber\\
&\le\|f\|_\8\left\{\left(\left(\E|R_{1,T}|^{\ff {p_0} {p_0-1}}\right)^{\ff {p_0-1} {p_0}}+\left(\E |R_{2,T}|^{\ff {p_0} {p_0-1}}\right)^{\ff {p_0-1} {p_0}}\right)\right.\nonumber\\
&\qquad \times \Big(\E\Big|\int_0^T\<\si^{-1}(b(Y_s)-b(Y_{s_\de})),\d W_s\>\Big|^{p_0 }\Big)^{\ff 1 {p_0}}\nonumber\\
&\qquad +\ff{1}{2}\left(\left(\E|R_{1,T}|^{\ff {p_0+1} {p_0-1}}\right)^{\ff {p_0-1} {p_0+1}}+\left(\E |R_{2,T}|^{\ff {p_0+1} {p_0-1}}\right)^{\ff {p_0-1} {p_0+1}}\right)\nonumber\\
&\qquad\times \left.\int_0^T\Big(\E\big||\si^{-1}b(Y_{s_\de})|^2-|\si^{-1}b(Y_{s})|^2\big|^{\ff {p_0+1} 2}\Big)^{\ff 2 {p_0+1}}\d s\right\}.
\end{align}

Define the stopping time $\tau_{1,n}=\inf\{t>0: \int_0^t|\si^{-1}b(Y_s)|^2\d s\ge n\}$.
Then $\hat{W}_{n,t}=W_t-\int_0^{t\wedge\tau_{1,n}}|\si^{-1}b(Y_s)|^2\d s$ is a Brownian motion under $\Q_1$. By H\"older's inequality, we arrive at
\beg{align*}
\E R_{1,T\we \tau_{1,n}}^{\ff {p_0} {p_0-1}}&=\E\exp\left\{\ff {p_0} {p_0-1}\int_0^{T\we \tau_{1,n}}\<\si^{-1}b(Y_s),\d W_s\>\right.\\
&\qquad\qquad \qquad \left.-\ff {p_0} {2(p_0-1)}\int_0^{T\we \tau_{1,n}}|\si^{-1}b(Y_s)|^2\d s\right\}\\
&\le \left(\E M_{1, T\we \tau_{1,n}}\right)^{1/2}\left(\E\exp\left\{\ff {p_0(p_0+1)} {(p_0-1)^2}\int_0^{T\we \tau_{1,n}}|\si^{-1}b(Y_s)|^2\d s\right\}\right)^{1/2}\\
&=\left(\E\exp\left\{\ff {p_0(p_0+1)} {(p_0-1)^2}\int_0^{T\we \tau_{1,n}}|\si^{-1}b(Y_s)|^2\d s\right\}\right)^{1/2},
\end{align*}
where
\beg{align*}
M_{1,{T\we \tau_{1,n}}}& =\exp\left\{\ff {2p_0} {p_0-1}\int_0^{T\we \tau_{1,n}}\<\si^{-1}b(Y_s),\d W_s\>\right.\\
&\qquad\qquad \qquad \left.-2\left(\ff {p_0} {p_0-1}\right)^2\int_0^{T\we \tau_{1,n}}|\si^{-1}b(Y_s)|^2\d s\right\},
\end{align*}
which is an exponential martingale. Similarly, 
\beg{align*}
\E R_{1,T\we \ta_{1,n}}^{\ff {p_0+1} {p_0-1}}\leq \left(\E\exp\left\{\ff {(p_0+3)(p_0+1)} {(p_0-1)^2}\int_0^{T\we \tau_{1,n}}|\si^{-1}b(Y_s)|^2\d s\right\}\right)^{1/2}. 
\end{align*}
It is clear that $\ff {(p_0+3)(p_0+1)} {(p_0-1)^2}> \ff {p_0(p_0+1)} {(p_0-1)^2}$. Then, due to \eqref{ine-Tp0}, Lemma \ref{lemn} and \eqref{ine-exp}, we have by  letting $n\ra+\infty$ that
$$\E\left( R_{1,T }^{\ff {p_0+1} {p_0-1}}+R_{1,T}^{\ff {p_0} {p_0-1}}\right)<\infty.$$ 
Similarly, we can prove by Lemma \ref{lemm} that  
$$\E\left( R_{2,T }^{\ff {p_0+1} {p_0-1}}+R_{2,T}^{\ff {p_0} {p_0-1}}\right)<\infty.$$

By (H1), for $s\geq\de$, we have
\beg{align}\label{bY}
&\E|b(Y_s)-b(Y_{s_\de})|^{p_0}=\E|b(x+\si W_s)-b(x+\si W_{s_\de})|^{p_0}\nonumber\\
 &=\int_{\R^d}\int_{\R^d}\ff {|b(x+y)-b(x+z)|^{p_0}e^{-\ff {|\si^{-1}z|^2} {2s_\de}}e^{-\ff {|\si^{-1}(y-z)|^2} {2(s-s_\de)}}} {\sq{(2\pi s_\de )^d\det(\si\si^*)}\cdot \sq{(2\pi(s-s_\de))^d\det(\si\si^*)}} \d z\d y\nonumber\\
 &\le  \ff{\|\si\|^{2d}}{\pi^d\det(\si\si^*)}\int_{\R^d}\int_{\R^d}
 \ff{|b(v)-b(u)|^{p_0} e^{-\ff {|u-x|^2} {2s_\de\|\si\|^2}}e^{-\ff {|v-u|^2} {2(s-s_\de)\|\si\|^2}}} 
 {(2 s_\de \|\si\|^2)^{d/2}\cdot (2(s-s_\de)\|\si\|^2)^{d/2}} \d u\d v
\nonumber\\
 &\le \ff{\|\si\|^{2d}}{\pi^d\det(\si\si^*)}(\phi(2s_\de\|\si\|^2)(2(s-s_\de)\|\si\|^2)^\al)^{p_0}.
\end{align}
Noting that
\beg{align}\label{phi}
\lim_{\delta\rightarrow0+}\sum_{k=1}^{[T/\de] }\phi^2(2k\de\|\si\|^2)\de &=\int_0^T\phi^2(2\|\si\|^2r)\d r\nonumber\\
&=\ff{1}{2\|\si\|^2}\int_0^{2\|\si\|^2T}\phi^2(s)\d s<\8,
\end{align}
the BDG inequality, \eqref{bY}, \eqref{phi} and $p_0\geq 2$ imply that
\beg{align}\label{G1}
G_{1,T}
& =\Big(\E\Big|\int_0^T\<\si^{-1}(b(Y_s)-b(Y_{s_\de})),\d W_s\>\Big|^{p_0}\Big)^{1/p_0}\nonumber\\
&\le \left(\ff {p_0} {p_0-1}\right)^{\ff {p_0} 2}\left(\ff {p_0(p_0-1)} 2\right)^{\ff 1 2}\|\si^{-1}\|
\left(\int_0^T\big(\E|b(Y_s)-b(Y_{s_\de})|^{p_0 }\big)^{\ff 2 {p_0 }}\d s\right)^{\ff 1 2}\nonumber\\
&\le \de^\al\ff {2^{\ff {\al-1} 2}\|\si\|^{\ff {2d} {p_0}+\al-1}\|\si^{-1}\|} {\left(\pi^d\det(\si\si^*)\right)^{\ff 1 {p_0}}}\left(\ff {p_0} {p_0-1}\right)^{\ff {p_0} 2}\left(\ff {p_0(p_0-1)} 2\right)^{\ff 1 2}\int_0^{2\|\si\|^2T}\phi^2(s)\d s\nonumber\\
&= C_{T,p_0,\si,\al,\ph}\de^{\alpha}.
 \end{align}

Noting that for any $p\geq 1$
\beg{align}\label{eq-Ym}
\E|Y_t|^p\le2^{p-1}\left(|x|^p+(\sq t \|\si\|)^p\E| W_1|^p\right),
\end{align}
we have
\begin{align*}
&\left(\E|b(Y_s)+b(Y_{s_\de})|^{\ff {p_0(p_0+1)} {p_0-1}}\right)^{\ff {p_0-1} {p_0(p_0+1)}}\\
&\qquad \leq \left(\E\left(2L_1+L_2(|Y_s|+|Y_{s_\de}|)\right)^{\ff {p_0(p_0+1)} {p_0-1}}\right)^{\ff {p_0-1} {p_0(p_0+1)}}\\
&\qquad \leq 6\left\{L_1+2L_2\left(|x|+\sq {T}\|\si\|\left(\E |W_1|^{\ff {p_0(p_0+1)} {p_0-1}}\right)^{\ff {p_0-1} {p_0(p_0+1)}}\right)\right\}\\
&\qquad =: C_{T,p_0,\si,L_1,L_2,x}.
\end{align*}
Combining this with (H1), (H2) and taking the similar argument as in \eqref{G1}, we obtain   
 \beg{align}\label{G2}
 &G_{2,T}=\ff{1}{2}\int_0^T\left(\E\Big||\si^{-1}b(Y_{s_\de})|^2-|\si^{-1}b(Y_{s})|^2\Big|^{\ff {p_0+1} 2}\right)^{\ff 2 {p_0+1}}\d s\nonumber\\
 &\le\ff{\|\si^{-1}\|^2}{2}\int_0^T\Big(\E |b(Y_s)-b(Y_{s_\de})|^{\ff {p_0+1}  2}|b(Y_s)+b(Y_{s_\de})|^{\ff {p_0+1} 2}\Big)^{\ff 2 {p_0+1}}\d s\nonumber\\
 &\leq \ff{\|\si^{-1}\|^2}{2}\int_0^T\left(\E |b(Y_s)-b(Y_{s_\de})|^{p_0}\right)^{\ff 1 {p_0}}\left(\E|b(Y_s)+b(Y_{s_\de})|^{\ff {p_0(p_0+1)} {p_0-1}}\right)^{\ff {p_0-1} {p_0(p_0+1)}}\d s\nonumber\\
&\leq  \ff{\|\si^{-1}\|^2}{2}C_{T,p_0,\si,L_1,L_2,x}\int_0^T\left(\E |b(Y_s)-b(Y_{s_\de})|^{p_0}\right)^{\ff 1 {p_0}}\d s\nonumber\\
&\leq C_{T,p_0,\si,L_1,L_2,\ph,x}\de^\al,
\end{align} 
where 
\beg{align*}
C_{T,p_0,\si,L_1,L_2,\ph,x}=\ff {2^{\al-2}\|\si\|^{\ff {2d} {p_0}+2\al-2}\|\si^{-1}\|^2C_{T,p_0,\si,L_1,L_2,x} } {(\pi^d \det(\si\si^*))^{\ff 1 {p_0}}} \int_0^{2\|\si\|^2T}\ph(s)\d s.
\end{align*}
Therefore, the conclusion holds under (H1) and (H2).

Step (ii), we claim that if $b$ satisfies the sublinear condition, then convergence holds for any $T > 0$.  In fact, for any given $T>0$, we can always choose $L_2=\ep>0$ small enough such that \eqref{ine-Tp0} holds. Then  \eqref{eq2.3} holds.

The proof is therefore complete.
\end{proof}

\noindent\textbf{Acknowledgements}

\medskip

The third author was supported by the National Natural Science Foundation of China (Grant No. 11901604, 11771326).


{\footnotesize 
 \begin{thebibliography}{99}
 
 
\bibitem{BT}
V. Bally, D. Talay, The law of the Euler scheme for stochastic differential equations, I. Convergence rate of the distribution function, Probab. Theory Related Fields 104 (1) (1996) 43--60.
 
\bibitem{BHY} J. Bao, X. Huang, C.Yuan, Convergence rate of Euler-Maruyama scheme for SDEs with H\"older-Dini continuous drifts. J. Theoret. Probab. 32 (2019),  848--871. 

\bibitem{BHZ} J. Bao, X. Huang, S.-Q. Zhang, Convergence rate of EM algorithm for SDEs under integrability condition, 

\bibitem{B} J. Bao, J. Shao, Weak convergence of path-dependent SDEs with irregular coefficients, Arxiv:1809.03088


\bibitem{NP}N. Halidias, P.E. Kloeden, A note on the Euler-Maruyama scheme for stochastic differential equations with a discontinuous monotone drift coefficient, BIT Numer Math, 48 (2008),51--59.


\bibitem{KLY} A. Kohatsu-Higa, A. Lejay, K. Yasuda, Weak rate of convergence of the Euler-Maruyama scheme for stochastic differential equations with non-regular drift. J. Comp. Appl. Math. 326 (2017), 138--158.


\bibitem{SV}V. Konakov, S. Menozzi, Weak error for the Euler scheme approximation of diffusions with non-smooth coefficients, Elect. J. Probab. 22 (2017),1--47.

\bibitem{NM}N.V. Krylov, M. R\"ockner, Strong solutions of stochastic equations with singular time dependent drift, Probab. Theory Related Fields. 
131 (2005), 154--196.

\bibitem{GM1} G. Leobacher, M. Sz\"olgyenyi, A numerical method for SDEs with discontinuous drift, BIT Numer Math, 56(2016),151--162.

\bibitem{GM2} G. Leobacher, M. Sz\"olgyenyi, A strong order $1/2$ method for multidimensional SDEs with discontinuous drift, Ann. Appl. Probab., 27 (2017), 2383--2418.

\bibitem{LM}
	V. Lemaire and S. Menozzi, On Some non Asymptotic Bounds for the Euler Scheme, \textit{Electron. J. Probab.} {\bf 15} (2010), no 53, 1645--1681.

\bibitem{AML} A. Neuenkirch, M. Sz\"olgyenyi, L. Szpruch, An adaptive Euler-Maruyama scheme for stochastic differential equations with discontinuous drift and its convergence analysis, SIAM J. Numer. Anal., 57 (2019),378--403.

\bibitem{NT1} H. L. Ngo, D. Taguchi, Strong rate of convergence for Euler-Maruyama approximation of stochastic differential equations with irregular coefficients, Math. Comp., 85(2016),1793--1819.

\bibitem{NT2} H. L. Ngo, D. Taguchi, On the Euler-Maruyama approximation for one-dimensional stochastic differential equations with irregular coefficients, IMA J. Numer. Anal., 37 (2017), 1864--1883.

\bibitem{S} J. Shao, Weak convergence of Euler-Maruyama's approximation for SDEs under integrability condition, https://arxiv.org/abs/1808.07250


	
\bibitem{W18}
	F.-Y. Wang, Estimates for invariant probability measures of degenerate SPDEs with singular and path-dependent drifts. 
	\textit{Proba. Theory $\&$ Related Fields}  (2018), 1--34.

\bibitem{XZ16}
	L. Xie, X. Zhang, Sobolev differentiable flows of SDEs with local Sobolev and super-linear growth coefficients, {\it Ann. Probab.} {\bf 44}  (2016)  3661--3687.

\bibitem{XZ} L. Xie, X. Zhang, Ergodicity of stochastic differential equations with jumps and singular coefficients, https://arxiv.org/abs/1705.07402 
 
 


\bibitem{X1}X.  Zhang, Strong solutions of SDEs with singular drift and Sobolev diffusion coefficients, Stochastic Process. Appl. 115 (2015),1805--1818.



















\end{thebibliography}
}
\end{document}